\documentclass[a4paper,11pt]{amsart}
\usepackage{graphicx}
\usepackage{mathptmx}
\usepackage{amsmath}
\usepackage{amssymb}
\usepackage{enumitem}
\usepackage{xcolor}
\usepackage{xparse}
\NewDocumentCommand{\eulerian}{omm}
 {%
  \genfrac<>{0pt}{}{#2}{#3}%
  \IfValueT{#1}{_{\!#1}}%
 }

 \newmuskip\pFqmuskip

\newcommand*\pFq[6][8]{%
  \begingroup 
  \pFqmuskip=#1mu\relax
  \mathchardef\normalcomma=\mathcode`,
  \mathcode`\,=\string"8000
  \begingroup\lccode`\~=`\,
  \lowercase{\endgroup\let~}\pFqcomma
  {}_{#2}F_{#3}{\left(\genfrac..{0pt}{}{#4}{#5}\bigg|#6\right)}%
  \endgroup
}
\newcommand{\pFqcomma}{{\normalcomma}\mskip\pFqmuskip}

\newtheorem{theorem}{Theorem}
\newtheorem{lemma}[theorem]{Lemma}

\begin{document}

\title[Some identities on degenerate hyperharmonic numbers]{Some identities on degenerate hyperharmonic numbers}

\author{Taekyun  Kim}
\address{Department of Mathematics, Kwangwoon University, Seoul 139-701, Republic of Korea}
\email{tkkim@kw.ac.kr}

\author{DAE SAN KIM}
\address{Department of Mathematics, Sogang University, Seoul 121-742, Republic of Korea}
\email{dskim@sogang.ac.kr}


\maketitle

\begin{abstract}
The aim of this paper is to investigate some properties, recurrence relations and identities involving degenerate hyperharmonic numbers, hyperharmonic numbers and degenerate harmonic numbers. In particular, we derive an explicit expression of the degenerate hyperharmonic numbers in terms of the degenerate harmonic numbers. This is a degenerate version of the corresponding identity representing the hyperharmonic numbers in terms of harmonic numbers due to Conway and Guy.
\end{abstract}

\section{Introduction}
Recent explorations for various degenerate versions of quite a few special numbers and polynomials have been fascinating and fruitful, which began with the pioneering work of Carlitz in [2,3].
These quests for degenerat versons are not only limited to special polynomials and numbers but also extended to some transcendental functions, like gamma functions (see [11]). Many different tools are used, which include generating functions, combinatorial methods, $p$-adic analysis, umbral calculus, operator theory, differential equations, special functions, probability theory and analytic number theory (see [8-10,12,14,16] and the references therein). It is also worth mentioning that $\lambda$-umbral calculus has been introduced in [8], which turns out to be more convenient than the `classical' umbral calculus when dealing with degenerate Sheffer polynomials.\par
The aim of this paper is to investigate some properties, recurrence relations and identities involving degenerate hyperharmonic numbers, hyperharmonic numbers and degenerate harmonic numbers (see \eqref{6}, \eqref{9}, \eqref{12}). The novelty of this paper is the derivation of an explicit expression of the degenerate hyperharmonic numbers in terms of the degenerate harmonic numbers (see Theorem 3). This is a degenerate version of the corresponding identity representing the hyperharmonic numbers in terms of harmonic numbers due to Conway and Guy (see \eqref{10}).\par
The outline of this paper is as follows. In Section 1, we recall the degenerate exponentials and the degenerate logarithms. We remind the reader of the harmonic numbers and their generating function, and of the degenerate harmonic numbers and their generating function. Then we recall the hyperharmonic numbers due to Conway and Guy [5], its explicit expression in terms of harmonic numbers and their generating function. Finally, we remind the reader of the recently introduced degenerate hyperharmonic numbers, which are a degenerate version of the hyperharmonic numbers, and of their generating function.
Section 2 is the main result of this paper. We obtain an identity involving the degenerate hyperharmonic numbers and the hyperharmonic numbers in Theorem 2. It is obtained by taking higher order derivatives of the generating function of the degenerate hyperharmonic numbers in \eqref{13}. Theorem 3 is a degenerate version of the explicit expression for the hyperharmonic numbers \eqref{10}, which is obtained from Theorem 2 and Lemma 1 about explicit expressions of certain polynomials. In Section 3, we get an identity involving the degenerate hyperharmonic numbers and the degenerate zeta function. \par
For any nonzero $\lambda\in\mathbb{R}$, the degenerate exponential functions are defined by
\begin{equation}
e_{\lambda}^{x}(t)=(1+\lambda t)^{\frac{x}{\lambda}}=\sum_{n=0}^{\infty}\frac{(x)_{n,\lambda}}{n!}t^{n},\quad (\mathrm{see}\ [2,8,9,10,13]), \label{1}
\end{equation}
where
\begin{equation}
(x)_{0,\lambda}=1,\quad (x)_{n,\lambda}=x(x-\lambda)\cdots(x-(n-1)\lambda),\quad (n\ge 1),\quad(\mathrm{see}\ [8,9,10]).\label{2}	
\end{equation}
When $x=1$, we use the notation $e_{\lambda}(t)=e_{\lambda}^{1}(t)$. \par
Let $\log_{\lambda}t$ be the compositional inverse function of $e_{\lambda}(t)$ such that $\log_{\lambda}(e_{\lambda}(t))=e_{\lambda}\big(\log_{\lambda}(t)\big)=t$.
It is called the degenerate logarithm and given by
\begin{equation}
\log_{\lambda}(1+t)=\sum_{k=1}^{\infty}\frac{\lambda^{k-1}(1)_{k,\frac{1}{\lambda}}}{k!}t^{k}=\frac{1}{\lambda}\big((1+t)^{\lambda}-1\big),\quad (\mathrm{see}\ [6,7,9]). \label{3}	
\end{equation}
Note that $\displaystyle \lim_{\lambda\rightarrow 0}\log_{\lambda}(1+t)=\log(1+t)\displaystyle$ and $\displaystyle\lim_{\lambda\rightarrow 0}e_{\lambda}(t)=e^{t}\displaystyle$. \par
It is well known that the harmonic numbers are defined by
\begin{equation}
H_{0}=0,\quad H_{n}=1+\frac{1}{2}+\frac{1}{3}+\cdots+\frac{1}{n},\quad(n \ge 1),\quad (\mathrm{see}\ [3,4,5,11,12]).\label{4}
\end{equation}
From \eqref{4}, we can derive the generating function of harmonic numbers given by
\begin{equation}
-\frac{1}{1-t}\log(1-t)=\sum_{n=1}^{\infty}H_{n}t^{n},\quad (\mathrm{see}\ [5]). \label{5}	
\end{equation}
Recently, the degenerate harmonic numbers are defined by
\begin{equation}
H_{0,\lambda}=0,\quad H_{n,\lambda}=\sum_{k=1}^{n}\frac{1}{\lambda}\binom{\lambda}{k}(-1)^{k-1},\quad (n\in\mathbb{N}),\quad (\mathrm{see}\ [6]). \label{6}
\end{equation}
Note that $\displaystyle\lim_{\lambda\rightarrow 0}H_{n,\lambda}=H_{n}\displaystyle$. \par
From \eqref{3} and \eqref{6}, we can derive the generating function of the degenerate harmonic numbers given by
\begin{equation}
-\frac{1}{1-t}\log_{\lambda}(1-t)=\sum_{n=1}^{\infty}H_{n,\lambda}t^{n},\quad (\mathrm{see}\ [6]).\label{7}
\end{equation}
In 1996, Conway and Guy introduced the hyperharmoinc numbers $H_{n}^{(r)},\,(n,r \ge 0),$ which are defined by
\begin{equation}
H_{0}^{(r)}=0;\,\, H_{n}^{(0)}=\frac{1}{n},\,\, H_{n}^{(1)}=H_{n},\,\, H_{n}^{(r)}=\sum_{k=1}^{n}H_{k}^{(r-1)},\,\,(r\ge 2),\quad (\mathrm{see}\ [3,5,11,14,16]). \label{9}
\end{equation}
By \eqref{9}, we get
\begin{equation}
H_{n}^{(r)}=\binom{n+r-1}{r-1}(H_{n+r-1}-H_{r-1}),\quad (r\ge 1),\label{10}
\end{equation}
\begin{equation}
-\frac{\log(1-t)}{(1-t)^{r}}=\sum_{n=1}^{\infty}H_{n}^{(r)}t^{n},\quad (\mathrm{see}\ [5]).\label{11}
\end{equation}
In [12], Kim-Kim introduced the degenerate hyperharmonic numbers $H_{n,\lambda}^{(r)},\,\,(n \ge 0, r \ge 1),$ which are given by
\begin{equation}
H_{0,\lambda}^{(r)}=0;\quad H_{n,\lambda}^{(1)}=H_{n,\lambda},\quad H_{n,\lambda}^{(r)}=\sum_{k=1}^{n}	H_{k,\lambda}^{(r-1)},\quad (r\ge 2). \label{12}
\end{equation}
From \eqref{12}, we note that
\begin{align}
-\frac{\log_{\lambda}(1-t)}{(1-t)^{r}}&=\frac{1}{(1-t)^{r-1}}\bigg(\sum_{n=1}^{\infty}H_{n,\lambda}t^{n}\bigg)=\frac{1}{(1-t)^{r-2}}\bigg(\sum_{n=1}^{\infty}\bigg(\sum_{k=1}^{n}H_{k,\lambda}\bigg)\bigg)t^{n}\label{13} \\
&=\frac{1}{(1-t)^{r-2}}\sum_{n=1}^{\infty}H_{n,\lambda}^{(2)}t^{n}=\cdots=\sum_{n=1}^{\infty}H_{n,\lambda}^{(r)}t^{n},\quad (r\ge 1). \nonumber	
\end{align}
By \eqref{13}, we easily get
\begin{equation}
H_{k,\lambda}^{(0)}=\frac{1}{k!}\lambda^{k-1}(-1)^{k-1}(1)_{k,\frac{1}{\lambda}},\quad (k\ge 1).
\end{equation}

\section{Some identities on degenerate hyperharmonic numbers}
First, let us define the polynomial $q_{n}(\lambda)\in Q[\lambda],\ (n\ge 0)$, with $\deg q_{n}(\lambda)=n-1$ as
\begin{equation}
\bigg(1-\frac{\lambda}{r}\bigg)\bigg(1-\frac{\lambda}{r+1}\bigg)\cdots\bigg(1-\frac{\lambda}{r+n-1}\bigg)=1+\lambda q_{n}(\lambda),\quad (n\in\mathbb{N}, r\in\mathbb{N}). \label{14}	
\end{equation}
From \eqref{14}, we have
\begin{align}
& \bigg(1-\frac{\lambda}{r}\bigg)\bigg(1-\frac{\lambda}{r+1}\bigg)\cdots\bigg(1-\frac{\lambda}{r+n-1}\bigg)\label{15}\\
&=\frac{(-1)^{n}}{r(r+1)\cdots(r+n-1)}(\lambda-r)(\lambda-r-1)(\lambda-r-2)\cdots(\lambda-r-(n-1))\nonumber\\
&=\frac{(r-1)!(-1)^{n}}{(r+n-1)!}(\lambda-r)_{n}=\frac{(-1)^{n}}{\binom{r+n-1}{n}n!}(\lambda-r)_{n}=\frac{(-1)^{n}}{\binom{r+n-1}{n}}\binom{\lambda-r}{n},\nonumber	
\end{align}
where
\begin{displaymath}
	(x)_{0}=1,\quad (x)_{n}=x(x-1)\cdots(x-n+1),\quad (n\ge 1).
\end{displaymath}
For $n\ge 0$, the Stirling numbers of the first kind are defined by
\begin{equation}
(x)_{n}	=\sum_{k=0}^{n}S_{1}(n,k)x^{k},\quad (\mathrm{see}\ [1-16]). \label{16}
\end{equation}
By \eqref{16}, we get
\begin{align}
&\frac{(-1)^{n}}{r(r+1)\cdots(r+n-1)}(\lambda-r)(\lambda-r-1)\cdots(\lambda-r-n+1)=\frac{1}{(-r)_{n}}\sum_{l=0}^{n}(\lambda-r)^{l}S_{1}(n,l) \label{17}\\
&=\frac{1}{(-r)_{n}}\sum_{l=0}^{n}S_{1}(n,l)\sum_{k=0}^{l}\binom{l}{k}(-r)^{l-k}\lambda^{k}=\frac{1}{(-r)_{n}}\sum_{k=0}^{n}\lambda^{k}\sum_{l=k}^{n}S_{1}(n,l)\binom{l}{k}(-1)^{l-k}r^{l-k} \nonumber \\
&=1+\frac{1}{(-r)_{n}}\sum_{k=1}^{n}\lambda^{k}\sum_{l=k}^{n}\binom{l}{k}S_{1}(n,l)(-1)^{l-k}r^{l-k} \nonumber \\
&=1+\lambda\bigg(\frac{1}{(-r)_{n}}\sum_{k=1}^{n}\lambda^{k-1}\sum_{l=k}^{n}\binom{l}{k}S_{1}(n,l)(-1)^{l-k}r^{l-k}\bigg)\nonumber \\
&=1+\lambda q_{n}(\lambda). \nonumber
\end{align}
From \eqref{15} and \eqref{17}, we note that
\begin{equation}
1+\lambda q_{n}(\lambda)=\frac{(-1)^{n}}{\binom{r+n-1}{n}}\binom{\lambda-r}{n},\quad (n\ge 1,\ r\ge 1). \label{18}
\end{equation}
Therefore, by \eqref{18}, we obtain the following lemma.
\begin{lemma}
	For $r,n\in\mathbb{N}$, the polynomials $q_{n}(\lambda)\in Q[\lambda]$ with $\deg q_{n}(\lambda)=n-1$ are defined by
	\begin{displaymath}
		\prod_{i=0}^{n-1}\bigg(1-\frac{\lambda}{r+i}\bigg)=1+\lambda q_{n}(\lambda).
	\end{displaymath}
	Then we have
	\begin{displaymath}
		q_{n}(\lambda)=\frac{1}{\lambda}\bigg\{\frac{(-1)^{n}}{\binom{r+n-1}{n}}\binom{\lambda-r}{n}-1\bigg\}.
	\end{displaymath}
\end{lemma}
From \eqref{13}, we note that
\begin{align}
\frac{d^{k}}{dz^{k}}\bigg(-\frac{\log_{\lambda}(1-z)}{(1-z)^{r}}\bigg)&=\frac{d^{k}}{dz^{k}}\sum_{n=0}^{\infty}H_{n,\lambda}^{(r)}z^{n} \label{19}\\
&=\sum_{n=k}^{\infty}H_{n,\lambda}^{(r)}n(n-1)\cdots(n-k+1)z^{n-k} \nonumber \\
&=k!\sum_{n=0}^{\infty}H_{n+k,\lambda}^{(r)}\binom{n+k}{k}z^{n},\nonumber	
\end{align}
where $k,r$ are positive integers. \par
Now, we observe that
\begin{align}
\frac{d}{dz}\bigg(-\frac{\log_{\lambda}(1-z)}{(1-z)^{r}}\bigg)&=\frac{(1-z)^{\lambda}}{(1-z)^{r+1}}-\frac{r}{(1-z)^{r+1}}\log_{\lambda}(1-z) \label{20}\\
&=\frac{r}{(1-z)^{r+1}}\bigg\{\frac{\lambda}{r}\log_{\lambda}(1-z)+\frac{1}{r}-\log_{\lambda}(1-z)\bigg\}\nonumber \\
&=\frac{r}{(1-z)^{r+1}}\bigg\{\frac{1}{r}-\log_{\lambda}(1-z)\bigg(1-\frac{\lambda}{r}\bigg)\bigg\}.\nonumber
\end{align}
Thus, by \eqref{20}, we get
\begin{align}
\frac{d}{dz}\bigg(-\frac{\log_{\lambda}(1-z)}{(1-z)^{r}}\bigg)&=\frac{r}{(1-z)^{r+1}}\bigg\{\frac{1}{r}-\log_{\lambda}(1-z)\bigg(1-\frac{r}{\lambda}\bigg)\bigg\}.\label{21}
\end{align}
From \eqref{21}, we have
\begin{align}
&\frac{d^{2}}{dz^{2}}\bigg(-\frac{\log_{\lambda}(1-z)}{(1-z)^{r}}\bigg)=\frac{d}{dz}\bigg(\frac{r}{(1-z)^{r+1}}\bigg(\frac{1}{r}-\bigg(1-\frac{\lambda}{r}\bigg)\log_{\lambda}(1-z)\bigg)\bigg) \label{22} \\
&=\frac{r(r+1)}{(1-z)^{r+2}}\bigg(\frac{1}{r}-\bigg(1-\frac{\lambda}{r}\bigg)\log_{\lambda}(1-z)\bigg)+\frac{r}{(1-z)^{r+1}}\bigg(1-\frac{\lambda}{r}\bigg)\frac{(1-z)^{\lambda}}{1-z} \nonumber \\
&=\frac{r(r+1)}{(1-z)^{r+2}}\bigg(\frac{1}{r}-\log_{\lambda}(1-z)\bigg(1-\frac{\lambda}{r}\bigg)\bigg)+\frac{r(r+1)}{(1-z)^{r+2}}\bigg(1-\frac{\lambda}{r}\bigg)\bigg(\frac{\lambda}{r+1}\log_{\lambda}(1-z)+\frac{1}{r+1}\bigg) \nonumber \\
&=\frac{r(r+1)}{(1-z)^{r+2}}\bigg\{\frac{1}{r}+\frac{1}{r+1}\bigg(1-\frac{\lambda}{r}\bigg)-\bigg(1-\frac{\lambda}{r}\bigg)\bigg(1-\frac{\lambda}{r+1}\bigg)\log_{\lambda}(1-z)\bigg\}.\nonumber
\end{align}
Thus, by \eqref{22}, we get
\begin{equation}
\begin{aligned}
	&\frac{d^{2}}{dz^{2}}\bigg(-\frac{\log_{\lambda}(1-z)}{(1-z)^{r}}\bigg)\\
	&=\frac{r(r+1)}{(1-z)^{r+2}}\bigg\{\frac{1}{r}+\frac{1}{r+1}\bigg(1-\frac{\lambda}{r}\bigg)-\bigg(1-\frac{\lambda}{r}\bigg)\bigg(1-\frac{\lambda}{r+1}\bigg)\log_{\lambda}(1-z)\bigg\}.
\end{aligned}\label{23}
\end{equation}
From \eqref{23}, we note that
\begin{align}
&\frac{d^{3}}{dz^{3}}\bigg(-\frac{\log_{\lambda}(1-z)}{(1-z)^{r}}\bigg)=\frac{d}{dz}\bigg(\frac{d^{2}}{dz^{2}}\bigg(-\frac{\log_{\lambda}(1-z)}{(1-z)^{r}}\bigg)\bigg) \label{24}	\\
&=\frac{d}{dz}\bigg\{\frac{r(r+1)}{(1-z)^{r+2}}\bigg(\frac{1}{r}+\frac{1}{r+1}\bigg(1-\frac{\lambda}{r}\bigg)-\bigg(1-\frac{\lambda}{r}\bigg)\bigg(1-\frac{\lambda}{r+1}\bigg)\log_{\lambda}(1-z)\bigg\} \nonumber \\
&=\frac{r(r+1)(r+2)}{(1-z)^{r+3}}\bigg\{\frac{1}{r}+\frac{1}{r+1}\bigg(1-\frac{\lambda}{r}\bigg)-\bigg(1-\frac{\lambda}{r}\bigg)\bigg(1-\frac{\lambda}{r+1}\bigg)\log_{\lambda}(1-z)\bigg\}\nonumber \\
&\quad +\frac{r(r+1)}{(1-z)^{r+2}}\bigg(1-\frac{\lambda}{r}\bigg)\bigg(1-\frac{\lambda}{r+1}\bigg)\frac{\lambda}{1-z}\frac{1}{\lambda}\big((1-z)^{\lambda}-1+1\big) \nonumber \\
&=\frac{r(r+1)(r+2)}{(1-z)^{r+3}}\bigg\{\frac{1}{r}+\frac{1}{r+1}\bigg(1-\frac{\lambda}{r}\bigg)-\bigg(1-\frac{\lambda}{r}\bigg)\bigg(1-\frac{\lambda}{r+1}\bigg)\log_{\lambda}(1-z)\bigg\} \nonumber \\
&\quad +\frac{r(r+1)(r+2)}{(1-z)^{r+3}}\bigg(1-\frac{\lambda}{r}\bigg)\bigg(1-\frac{\lambda}{r+1}\bigg)\bigg(\frac{\lambda}{r+2}\log_{\lambda}(1-z)+\frac{1}{r+2}\bigg) \nonumber\\
&=\frac{r(r+1)(r+2)}{(1-z)^{r+3}}\bigg\{\frac{1}{r}+\frac{1}{r+1}\bigg(1-\frac{\lambda}{r}\bigg)+\frac{1}{r+2}\bigg(1-\frac{\lambda}{r}\bigg)\bigg(1-\frac{\lambda}{r+1}\bigg)\nonumber \\
&\quad -\bigg(1-\frac{\lambda}{r}\bigg)\bigg(1-\frac{\lambda}{r+1}\bigg)\bigg(1-\frac{\lambda}{r+2}\bigg)\log_{\lambda}(1-z)\bigg\}.\nonumber
\end{align}
Thus, by \eqref{24}, we get
\begin{equation}
\begin{aligned}
	&\frac{d^{3}}{dz^{3}}\bigg(-\frac{\log_{\lambda}(1-z)}{(1-z)^{r}}\bigg)\\
	&=\frac{r(r+1)(r+2)}{(1-z)^{r+3}}\bigg\{\frac{1}{r}+\frac{1}{r+1}\bigg(1-\frac{\lambda}{r}\bigg)+\frac{1}{r+2}\bigg(1-\frac{\lambda}{r}\bigg)\bigg(1-\frac{\lambda}{r+1}\bigg)\\
	&\quad -\bigg(1-\frac{\lambda}{r}\bigg)\bigg(1-\frac{\lambda}{r+1}\bigg)\bigg(1-\frac{\lambda}{r+2}\bigg)\log_{\lambda}(1-z)\bigg\}.
	\end{aligned}	
\end{equation}
Continuing this process, we have
\begin{align}
\frac{d^{k}}{dz^{k}}\bigg(-\frac{\log_{\lambda}(1-z)}{(1-z)^{r}}\bigg)&=\frac{r(r+1)\cdots(r+k-1)}{(1-z)^{r+k}}\bigg(\frac{1}{r}+\sum_{l=2}^{k}\frac{1}{r+l-1}\prod_{j=0}^{l-2}\bigg(1-\frac{\lambda}{r+j}\bigg)\bigg) \label{25} \\
&\quad -	r(r+1)\cdots(r+k-1)\prod_{l=0}^{k-1}\bigg(1-\frac{\lambda}{r+l}\bigg)\frac{\log_{\lambda}(1-z)}{(1-z)^{r+k}}. \nonumber
\end{align}
From Lemma 1 and \eqref{25}, we can derive the following equation. \eqref{26}
\begin{align}
&\frac{d^{k}}{dz^{k}}\bigg(-\frac{\log_{\lambda}(1-z)}{(1-z)r}\bigg) \label{26}	\\
&=\frac{r(r+1)\cdots(r+k-1)}{(1-z)^{r+k}}\bigg(\frac{1}{r}+\sum_{l=2}^{k}\frac{1}{r+l-1}\big(1+\lambda q_{l-1}(\lambda)\big)\bigg)\nonumber\\
&\quad -r(r+1)\cdots(r+k-1)(1+\lambda q_{k}(\lambda)\big)\frac{\log_{\lambda}(1-z)}{(1-z)^{k+r}} \nonumber \\
&=\frac{(r+k-1)!}{(r-1)!}\frac{1}{(1-z)^{r+k}}\bigg(\sum_{l=1}^{k}\frac{1}{r+l-1}+\lambda\sum_{l=2}^{k}\frac{q_{l-1}(\lambda)}{r+l-1}\bigg) \nonumber \\
&\quad -\frac{(r+k-1)!}{(r-1)!}\big(1+\lambda q_{k}(\lambda)\big)\frac{\log_{\lambda}(1-z)}{(1-z)^{k+r}} \nonumber \\
&=\frac{(r+k-1)!}{(r-1)!}\frac{1}{(1-z)^{r+k}}\bigg(\sum_{l=1}^{k}\frac{1}{r+l-1}-\log_{\lambda}(1-z)\bigg) \nonumber \\
&\quad +\lambda\frac{(r+k-1)!}{(r-1)!}\frac{1}{(1-z)^{r+k}}\bigg(\sum_{l=2}^{k}\frac{q_{l-1}(\lambda)}{r+l-1}-q_{k}(\lambda)\log_{\lambda}(1-z)\bigg) \nonumber \\
&=\frac{(r+k-1)!}{(r-1)!}\bigg(\sum_{n=0}^{\infty}\binom{n+r+k-1}{n}\bigg(H_{k+r-1}-H_{r-1}\bigg)z^{n}-\frac{\log_{\lambda}(1-z)}{(1-z)^{r+k}}\bigg)\nonumber\\
&\quad +\lambda\frac{(r+k-1)!}{(r-1)!}\bigg(\sum_{n=0}^{\infty}\binom{n+r+k-1}{n}\sum_{l=2}^{k}\frac{q_{l-1}(\lambda)}{r+l-1}z^{n}-q_{k}(\lambda)\frac{\log_{\lambda}(1-z)}{(1-z)^{r+k}}\bigg). \nonumber
\end{align}
From the generating function of degenerate hyperharmonic numbers in \eqref{13} and \eqref{26}, we obtain
\begin{align}
&\frac{d^{k}}{dz^{k}}\bigg(-\frac{\log_{\lambda}(1-z)}{(1-z)^{r}}\bigg) \label{27}	\\
&=\frac{(r+k-1)!}{(r-1)!}\bigg(\sum_{n=0}^{\infty}\binom{n+r+k-1}{n})(H_{k+r-1}-H_{r-1})z^{n}+\sum_{n=1}^{\infty}H_{n,\lambda}^{(k+r)}z^{n}\bigg)\nonumber \\
&\quad +\lambda\frac{(r+k-1)!}{r-1}\bigg(\sum_{n=0}^{\infty}\binom{n+r+k-1}{n}\sum_{l=2}^{k}\frac{q_{l-1}(\lambda)}{r+l-1}z^{n}+q_{k}(\lambda)\sum_{n=1}^{\infty}H_{n,\lambda}^{(r+k)}z^{n}\bigg) \nonumber \\
&=k!\sum_{n=0}^{\infty}\bigg\{\binom{n+r+k-1}{n}\binom{r+k-1}{k}(H_{k+r-1}-H_{r-1})+\binom{r+k-1}{k}H_{n,\lambda}^{(k+r)}\bigg\}z^{n}\nonumber\\
&\quad +\lambda k!\sum_{n=0}^{\infty}\bigg\{\binom{n+r+k-1}{n}\binom{r+k-1}{k}\sum_{l=2}^{k}\frac{q_{l-1}(\lambda)}{r+l-1}+q_{k}(\lambda)\binom{r+k-1}{k}H_{n,\lambda}^{(k+r)}\bigg\}z^{n}.\nonumber
\end{align}
From \eqref{10} and \eqref{27}, we note that
\begin{align}
&\frac{d^{k}}{dz^{k}}\bigg(-\frac{\log_{\lambda}(1-z)}{(1-z)^{r}}\bigg) \label{28} \\
&=k!\sum_{n=0}^{\infty}\bigg\{\binom{n+r+k-1}{n}H_{k}^{(r)}+\binom{r+k-1}{k}H_{n,\lambda}^{(k+r)}\bigg\}z^{n}\nonumber \\
&\quad +\lambda k!\sum_{n=0}^{\infty}\bigg\{\binom{n+r+k-1}{n}\binom{r+k-1}{k}\sum_{l=2}^{k}\frac{q_{l-1}(\lambda)}{r+l-1}+\binom{r+k-1}{k}q_{k}(\lambda)H_{n,\lambda}^{(k+r)}\bigg\}z^{n}\nonumber \\
&=k!\sum_{n=0}^{\infty}\bigg\{\binom{n+r+k-1}{n}H_{k}^{(r)}+\binom{r+k-1}{k}H_{n,\lambda}^{(k+r)} \nonumber \\
&\qquad +\lambda\bigg(\binom{n+r+k-1}{n}\binom{r+k-1}{k}\sum_{l=2}^{k}\frac{q_{l-1}(\lambda)}{r+l-1}+\binom{r+k-1}{k}q_{k}(\lambda)H_{n,\lambda}^{(k+r)}\bigg)\bigg\}z^{n}.\nonumber
\end{align}
Therefore, by \eqref{19} and \eqref{28}, we obtain the following theorem.
\begin{theorem}
For $r,k\in\mathbb{N}$, we have the following identity
\begin{align*}
&\binom{n+k}{k}H_{n+k,\lambda}^{(r)}=\binom{n+r+k-1}{n}H_{k}^{(r)}+\binom{r+k-1}{k}H_{n,\lambda}^{(k+r)} \quad\\
&+\lambda\bigg\{\binom{n+r+k-1}{n}\binom{r+k-1}{k}\sum_{l=2}^{k}\frac{q_{l-1}(\lambda)}{r+l-1}+\binom{r+k-1}{k}q_{k}(\lambda)H_{n,\lambda}^{(k+r)}\bigg\},
\end{align*}
where $q_{n}(\lambda)$ is the polynomial of degeree $n-1$ given by $q_{n}(\lambda)=\frac{1}{\lambda}\Big(\frac{(-1)^{n}}{\binom{n+r-1}{n}}\binom{\lambda-r}{n}-1\Big)$.
\end{theorem}
From Theorem 2 and Lemma 1, we have
\begin{align}
&\binom{n+k}{k}H_{n+k,\lambda}^{(r)}=\binom{n+r+k-1}{n}H_{k}^{(r)}+\binom{r+k-1}{k}H_{n,\lambda}^{(k+r)}\label{30}\\
&+\binom{n+r+k-1}{n}\binom{r+k-1}{k}\sum_{l=2}^{k}\frac{1}{r+l-1}\bigg(\frac{(-1)^{l-1}}{\binom{r+l-2}{l-1}}\binom{\lambda-r}{l-1}-1\bigg) \nonumber \\
&+\binom{r+k-1}{k}\bigg(\frac{(-1)^{k}}{\binom{r+k-1}{k}}\binom{\lambda-r}{k}-1\bigg)H_{n,\lambda}^{(k+r)}\nonumber \\
&=\binom{n+r+k-1}{n}H_{k}^{(r)}+\binom{n+r+k-1}{n}\binom{r+k-1}{k}
\nonumber\\
&\times \bigg(\sum_{l=2}^{k}\frac{(-1)^{l-1}\binom{\lambda-r}{l-1}}{\binom{r+l-2}{l-1}(r+l-1)}-\sum_{l=2}^{k}\frac{1}{r+l-1}\bigg)
+(-1)^{k}\binom{\lambda-r}{k}H_{n,\lambda}^{(k+r)}.\nonumber
\end{align}
Let us take $r=1$ in \eqref{30}. Then we have
\begin{align}
\binom{n+k}{k}H_{n+k,\lambda}&=\binom{n+k}{n}H_{k}+\binom{n+k}{n}\bigg(\sum_{l=2}^{k}\frac{1}{l}\binom{\lambda-1}{l-1}(-1)^{l-1}-\sum_{l=2}^{k}\frac{1}{l}\bigg) \label{31} \\
&\quad +(-1)^{k}\binom{\lambda-1}{k}H_{n,\lambda}^{(k+1)} \nonumber \\
&=\binom{n+k}{n}H_{k}+\binom{n+k}{k}\bigg(\sum_{l=1}^{k}\frac{1}{l}\binom{\lambda-1}{l-1}(-1)^{l-1}-\sum_{l=1}^{k}\frac{1}{l}\bigg) \nonumber \\
&\quad +(-1)^{k}\binom{\lambda-1}{k}H_{n,\lambda}^{(k+1)}\nonumber \\
&=\binom{n+k}{n}H_{k}+\binom{n+k}{k}\bigg(\sum_{l=1}^{k}\frac{1}{\lambda}\binom{\lambda}{l}(-1)^{l-1}-H_{k}\bigg)\nonumber\\
&\quad +(-1)^{k}\binom{\lambda-1}{k}H_{n,\lambda}^{(k+1)}
\nonumber\\
&=\binom{n+k}{n}H_{k}+\binom{n+k}{k}(H_{k,\lambda}-H_{k})+(-1)^{k}\binom{\lambda-1}{k}H_{n,\lambda}^{(k+1)} \nonumber \\
&=\binom{n+k}{n}H_{k,\lambda}+(-1)^{k}\binom{\lambda-1}{k}H_{n,\lambda}^{(k+1)}. \nonumber
\end{align}
Therefore, by \eqref{31}, we obtain the following theorem.
\begin{theorem}
For $n,k\in\mathbb{N}$, we have
\begin{displaymath}
H_{n,\lambda}^{(k+1)}=\frac{(-1)^{k}}{\binom{\lambda-1}{k}}\binom{n+k}{n}\big(H_{n+k,\lambda}-H_{k,\lambda}\big).
\end{displaymath}	
\end{theorem}
In Theorem 3, letting $\lambda \rightarrow 0$ gives the result in \eqref{10}. Namely, we have
\begin{displaymath}
H_{n}^{(k+1)}=\binom{n+k}{n}(H_{n+k}-H_{k}).
\end{displaymath}
This shows that Theorem 3 is a degenerate version of the expression in \eqref{10}.

\section{Further Remark}
For $\mathrm{Re}(\delta)>0$, the degenerate Hurwitz zeta function is defined by Kim-Kim as
\begin{displaymath}
\zeta_{\lambda}(s,\delta)=\sum_{n=0}^{\infty}\frac{(1)_{n,\lambda}}{(n+\delta)^{s}},\quad (\mathrm{Re}(s)>1),\quad (\mathrm{see}\ [7]).
\end{displaymath}
In particular, $\delta=1$, $\zeta_{\lambda}(s)=\zeta_{\lambda}(s,1)$ is called the degenerate zeta function. That is,
\begin{displaymath}
	\zeta_{\lambda}(s)=\sum_{n=1}^{\infty}\frac{(1)_{n-1,\lambda}}{n^{s}},\quad (\mathrm{Re}(s)>1),\quad (\mathrm{see}\ [7]).
\end{displaymath}
\begin{align}
&\frac{H_{1,\lambda}^{(r)}}{1^{m}}+\frac{H_{2,\lambda}^{(r)}}{2^{m}}(1)_{1,\lambda}+\frac{H_{3,\lambda}^{(r)}}{3^{m}}(1)_{2,\lambda}+\frac{H_{4,\lambda}^{(r)}}{4^{m}}(1)_{3,\lambda}+\cdots \label{32}\\
&=\frac{H_{1,\lambda}^{(r-1)}}{1^{m}}+\frac{H_{1,\lambda}^{(r-1)}+H_{2,\lambda}^{(r-1)}}{2^{m}}(1)_{1,\lambda}+\frac{1}{3^{m}}(H_{1,\lambda}^{(r-1)}+H_{2,\lambda}^{(r-1)}+H_{3,\lambda}^{(r-1)})(1)_{2,\lambda}+\cdots \nonumber \\
&=H_{1,\lambda}^{(r-1)}\sum_{k=1}^{\infty}\frac{(1)_{k-1,\lambda}}{k^{m}}+H_{2,\lambda}^{(r-1)}\sum_{k=1}^{\infty}\frac{(1)_{k,\lambda}}{(k+1)^{m}}+H_{3,\lambda}^{(r-1)}\sum_{k=1}^{\infty}\frac{(1)_{k+1,\lambda}}{(k+2)^{m+1}}+\cdots \nonumber  \\
&=\sum_{n=1}^{\infty}H_{n,\lambda}^{(r-1)}\sum_{k=1}^{\infty}\frac{(1)_{k+n-2,\lambda}}{(k+n-1)^{m}}=\sum_{n=1}^{\infty}H_{n,\lambda}^{(r-1)}\sum_{k=0}^{\infty}\frac{(1)_{k+n-1,\lambda}}{(k+n)^{m}}\nonumber \\
&=\sum_{n=1}^{\infty}H_{n,\lambda}^{(r-1)}\sum_{k=n}^{\infty}\frac{(1)_{k-1,\lambda}}{k^{m}}=-\sum_{n=1}^{\infty}H_{n,\lambda}^{(r-1)}\sum_{l=1}^{n-1}\frac{(1)_{l-1,\lambda}}{l^{m}}+\sum_{n=1}^{\infty}H_{n,\lambda}^{(r-1)}\sum_{k=1}^{\infty}\frac{(1)_{k-1,\lambda}}{k^{m}}\nonumber \\
&=-\sum_{n=1}^{\infty}H_{n,\lambda}^{(r-1)}\sum_{l=1}^{n-1}\frac{(1)_{l-1,\lambda}}{l^{m}}+\sum_{n=1}^{\infty}H_{n,\lambda}^{(r-1)}\zeta_{\lambda}(m). 	\nonumber
\end{align}
From \eqref{32}, we note that
\begin{displaymath}
\sum_{n=1}^{\infty}\bigg(\frac{H_{n,\lambda}^{(r)}}{n^{m}}(1)_{n-1,\lambda}+H_{n,\lambda}^{(r-1)}\sum_{l=1}^{n-1}\frac{(1)_{l-1,\lambda}}{l^{m}}\bigg)=\zeta_{\lambda}(m)\sum_{n=1}^{\infty}H_{n,\lambda}^{(r-1)}.
\end{displaymath}

\section{Conclusion}
In this paper, by using generating functions we investigated some properties, recurrence relations and identities involving degenerate hyperharmonic numbers, hyperharmonic numbers and degenerate harmonic numbers. The degenerate hyperharmonic numbers were introduced as a degenerate version of the hyperharmonic numbers which were introduced by Conway and Guy. In particular, derived was a degenerate version of the explicit expression of hyperharmonic numbers in terms of harmonic numbers, namely that of the degenerate hyperharmonic numbers in terms of the degenerate harmonic numbers.\par
It is one of our future projects to continue to study various degenerate versions of some special polynomials and numbers and to find their applications to physics, science and engineering as well as mathematics.

\end{document}